\documentclass[a4paper,11pt]{article}
\usepackage{amssymb}

\newcommand{\D}{\displaystyle}

\newcommand{\C}{\mathbb{C}}
\newcommand{\R}{\mathbb{R}}

\newcommand{\N}{\mathbb{N}}

\newcommand{\beq}{\begin{equation} }
\newcommand{\eqq}{\end{equation} }
\newcommand{\cuad}{{\sqcap\kern-.68em\sqcup}}

\newcommand{\norm}[1]{\|#1\|}

\newtheorem{definition}{Definition}[section]
\newtheorem{teo}{Theorem}[section]

\newtheorem{proposition}{Proposition}[section]

\newtheorem{lemma}{Lemma}[section]

\newtheorem{remark}{Remark}[section]
\newcommand{\bremark}{\begin{remark} \em}
\newcommand{\eremark}{\end{remark} }

\def\beeq{\begin{equation}}
\def\eeq{\end{equation}}
\newcommand{\begeqaet}{\begin{eqnarray*}}
\newcommand{\eneqaet}{\end{eqnarray*}}

\begin{document}

\begin{center}{\bf  \Large   Semilinear fractional elliptic equations with \qquad\medskip
 gradient nonlinearity involving measures\\
 }\medskip
\medskip
\bigskip
{\small
       {\bf Huyuan Chen}\footnote{chenhuyuan@yeah.net}\smallskip

       Departamento de Ingenier\'{\i}a  Matem\'atica
\\ Universidad de Chile, Chile\\[2mm]
     {\bf Laurent V\'{e}ron}\footnote{Laurent.Veron@lmpt.univ-tours.fr}

\smallskip
Laboratoire de Math\'{e}matiques et Physique Th\'{e}orique
\\  Universit\'{e} Fran\c{c}ois Rabelais, Tours, France}\\[1mm]
\bigskip

\medskip
\begin{abstract}
We study the existence of solutions to the fractional elliptic equation (E1) $(-\Delta)^\alpha
u+\epsilon g(|\nabla u|)=\nu $ in an open bounded regular domain $\Omega$ of $\R^N (N\ge2)$,
subject to the condition (E2) $u=0$ in $\Omega^c$, where $\epsilon=1$ or $-1$, $(-\Delta)^\alpha$
denotes the fractional Laplacian with $\alpha\in(1/2,1)$, $\nu$ is a
Radon measure and $g:\R_+\mapsto\R_+$ is a continuous  function. We prove the existence of weak solutions
for problem (E1)-(E2) when  $g$ is subcritical.
Furthermore, the asymptotic behavior and uniqueness of solutions are  described when $\epsilon=1$, $\nu$ is  Dirac mass and $g(s)=s^p$ with $p\in(0,\frac{N}{N-2\alpha+1})$.
\end{abstract}
\end{center}

\tableofcontents \vspace{1mm}
  \vspace{1mm}


\setcounter{equation}{0}
\section{Introduction}
Let $\Omega\subset \R^N(N\geq2)$ be an open bounded  $C^2$ domain and
$g:\R_+\mapsto\R_+$ be a continuous  function.  The purpose of this paper is to study the
existence of weak solutions  to the semilinear fractional elliptic problem with $\alpha\in(1/2,1)$,

\begin{equation}\label{eq1.1}
 \arraycolsep=1pt
\begin{array}{lll}
 (-\Delta)^\alpha  u+\epsilon g(|\nabla u|)=\nu\quad & \rm{in}\quad\Omega,\\[2mm]
 \phantom{   (-\Delta)^\alpha  u+\epsilon g(|\nabla u|}
u=0\quad & \rm{in}\quad \Omega^c,
\end{array}
\end{equation}
where $\epsilon=1$ or  $-1$  and $\nu\in\mathfrak{M}(\Omega,\rho^\beta)$ with $\beta\in[0,2\alpha-1)$. Here
$\rho(x)=dist(x,\Omega^c)$ and $\mathfrak{M}(\Omega,\rho^\beta)$ is the
space of Radon measures in $\Omega$ satisfying
\begin{equation}\label{radon}
\int_{\Omega}\rho^\beta d|\nu|<+\infty.
\end{equation}
In particular, we denote
$\mathfrak{M}^b(\Omega)=\mathfrak{M}(\Omega,\rho^0)$. The associated positive cones are respectively
$\mathfrak{M}_+(\Omega,\rho^\beta)$ and $\mathfrak{M}_+^b(\Omega)$. According to the value of $\epsilon$, we speak of an absorbing nonlinearity the case $\epsilon=1$ and a source nonlinearity the case $\epsilon=-1$. The operator
$(-\Delta)^\alpha $  is the fractional Laplacian defined as
$$(-\Delta)^\alpha  u(x)=\lim_{\varepsilon\to0^+} (-\Delta)_\varepsilon^\alpha u(x),$$
where for $\varepsilon>0$,
\begin{equation}\label{1.2}
(-\Delta)_\varepsilon^\alpha  u(x)=-\int_{\R^N}\frac{ u(z)-
u(x)}{|z-x|^{N+2\alpha}}\chi_\varepsilon(|x-z|) dz
\end{equation}
and
$$\chi_\varepsilon(t)=\left\{ \arraycolsep=1pt
\begin{array}{lll}
0,\quad & \rm{if}\quad t\in[0,\varepsilon],\\[2mm]
1,\quad & \rm{if}\quad t>\varepsilon.
\end{array}
\right.$$

In a pioneering work,  Brezis \cite{B12} (also see B\'{e}nilan and Brezis \cite{BB11}) studied the existence and uniqueness of the solution to
the semilinear Dirichlet elliptic problem
\begin{equation}\label{eq003}
 \arraycolsep=1pt
\begin{array}{lll}
 -\Delta  u+h(u)=\nu \quad & \rm{in}\quad\Omega,\\[2mm]
 \phantom{ -\Delta  +g(u)}
u=0  \quad & \rm{on}\quad \partial\Omega,
\end{array}
\end{equation}
where $\nu$ is a bounded
measure in $\Omega$ and the function $h$ is nondecreasing, positive on $(0,+\infty)$ and satisfies that
$$\int_1^{+\infty}(h(s)-h(-s))s^{-2\frac{N-1}{N-2}}ds<+\infty.$$
Later on, V\'{e}ron \cite{V} improved this result in replacing the Laplacian by more general uniformly elliptic second order differential operator, where $\nu\in\mathfrak{M}(\Omega,\rho^\beta)$ with $\beta\in[0,1]$ and
$h$ is a nondecreasing function satisfying
$$\int_1^{+\infty}(h(s)-h(-s))s^{-2\frac{N+\beta-1}{N+\beta-2}}ds<+\infty.$$
The general semilinear elliptic problems involving measures
 such as the equations involving boundary measures have been intensively studied;
it was initiated by Gmira and V\'{e}ron \cite{GV} and then
this subject has being extended in various ways, see  \cite{Hung,BV,MV1,MV2,MV3,MV4} for details and \cite{MV5} for a general panorama.
In a recent work, Nguyen-Phuoc and V\'{e}ron \cite{NV} obtained the existence of solutions to the viscous Hamilton-Jacobi  equation
\begin{equation}\label{eq1.111}
 \arraycolsep=1pt
\begin{array}{lll}
 -\Delta  u+h(|\nabla u|)=\nu\quad & \rm{in}\quad\Omega,\\[2mm]
 \phantom{   -\Delta  u+g(|\nabla u|}
u=0\quad & \rm{on}\quad \partial\Omega,
\end{array}
\end{equation}
when $\nu\in\mathfrak{M}^b(\Omega)$, $h$ is a continuous nondecreasing function vanishing at $0$ which satisfies
$$\int_1^{+\infty}h(s)s^{-\frac{2N-1}{N-1}}ds<+\infty.$$

 During the last years there has also been a renewed and increasing interest in the study of linear and nonlinear integro-differential
 operators, especially, the fractional Laplacian, motivated by great applications in physics and by important links on the theory
 of L\'{e}vy processes, refer to  \cite{CS1,CV1,CV2,CFQ,FQ2,RS,S1,S}. Many estimates of its Green kernel and generation formula  can be found in the references \cite{BKN,CS}.
Recently, Chen and  V\'{e}ron \cite{CV2}
studied the semilinear fractional elliptic equation
\begin{equation}\label{1.22}
 \arraycolsep=1pt
\begin{array}{lll}
 (-\Delta)^\alpha  u+h(u)=\nu\quad & \rm{in}\quad\Omega,\\[2mm]
 \phantom{   (-\Delta)^\alpha  +h(u)}
u=0\quad & \rm{in}\quad \Omega^c,
\end{array}
\end{equation}
where $\nu\in\mathfrak{M}(\Omega,\rho^\beta)$ with $ \beta\in[0,\alpha]$. We proved
  the existence and uniqueness of the
solution to (\ref{1.22}) when the function $h$ is nondecreasing and satisfies
$$\int_1^{+\infty}(h(s)-h(-s))s^{-1-k_{\alpha,\beta}}ds<+\infty,$$
where
\begin{equation}\label{eq 2.10}
k_{\alpha,\beta}=\left\{
\arraycolsep=1pt
\begin{array}{lll}
\frac{N}{N-2\alpha},\quad &\rm{if}\quad
\beta\in[0,\frac{N-2\alpha}N\alpha],\\[2mm]
\frac{N+\alpha}{N-2\alpha+\beta},\quad &\rm{if}\quad
\beta\in(\frac{N-2\alpha}N\alpha,\alpha].
\end{array}
\right.
\end{equation}

Our interest in this article is to investigate the existence of weak solutions to fractional equations involving nonlinearity in the gradient term and with Radon measure. In order the fractional Laplacian  be the dominant operator in terms of order of differentiation, it is natural to assume that $\alpha\in (1/2,1)$.
\begin{definition}\label{weak definition}
We say that $u$ is a weak solution of (\ref{eq1.1}), if $u\in
L^1(\Omega)$, $|\nabla u|\in L_{loc}^1(\Omega)$, $g(|\nabla u|)\in L^1(\Omega,\rho^\alpha dx)$  and
\begin{equation}\label{weak sense}
\int_\Omega [u(-\Delta)^\alpha\xi+\epsilon g(|\nabla u|)\xi]dx=\int_\Omega\xi
d\nu,\quad\  \forall\ \xi\in \mathbb{X}_{\alpha},
\end{equation}
where $\mathbb{X}_{\alpha}\subset C(\R^N)$ is the space of functions
$\xi$ satisfying:\smallskip

\noindent (i) $\rm{supp}(\xi)\subset\bar\Omega$,\smallskip

\noindent(ii) $(-\Delta)^\alpha\xi(x)$ exists for all $x\in \Omega$
and $|(-\Delta)^\alpha\xi(x)|\leq C$ for some $C>0$,\smallskip

\noindent(iii) there exist $\varphi\in L^1(\Omega,\rho^\alpha dx)$
and $\varepsilon_0>0$ such that $|(-\Delta)_\varepsilon^\alpha\xi|\le
\varphi$ a.e. in $\Omega$, for all
$\varepsilon\in(0,\varepsilon_0]$.\smallskip
\end{definition}

We denote by  $G_\alpha$ the Green kernel of $(-\Delta)^\alpha$ in
$\Omega\time\Omega $ and  by $\mathbb{G}_\alpha[.]$ the associated Green operator
defined by
\begin{equation}\label{optimal0}
\mathbb{G}_\alpha[\nu](x)=\int_{\Omega}G_\alpha(x,y) d\nu(y),\qquad\forall\ \nu\in
\mathfrak{M}(\Omega,\rho^\alpha).
\end{equation}
Using bounds of $\mathbb{G}_\alpha[\nu]$, we obtain in section 2  some crucial estimates which will play an important role in our construction of weak solutions. Our main result in the case $\epsilon=1$ is the following.

\begin{teo}\label{teo 1}
 Assume that $\epsilon=1$ and $g:\R_+\mapsto\R_+$
is a continuous function verifying $g(0)=0$ and
\begin{equation}\label{1.4}
\int_1^{+\infty} g(s)s^{-1-p^*_{\alpha}}ds<+\infty,
\end{equation}
where \begin{equation}\label{annex 00}
p^*_\alpha= \frac{N}{N-2\alpha+1}.
\end{equation}
Then for any $\nu\in\mathfrak{M}_+(\Omega,\rho^{\beta})$ with $\beta\in[0,2\alpha-1)$, problem
(\ref{eq1.1}) admits a nonnegative  weak solution $u_\nu$ which satisfies
\begin{equation}\label{1.5}
 u_\nu\le \mathbb{G}_\alpha[\nu].
\end{equation}
\end{teo}

As in the case $\alpha=1$, uniqueness remains an open question.
We remark  that the critical value $p^*_\alpha$ is independent of $\beta$.
A similar fact was first observed when dealing with problem (\ref{1.22}) where the critical value $k_{\alpha,\beta}$ defined by (\ref{eq 2.10}) does not depend on $\beta$ when $\beta\in[0,\frac{N-2\alpha}N\alpha]$.

\smallskip

When $\epsilon=-1$, we have to consider the critical value $ p^*_{\alpha,\beta}$ which depends truly on $\beta$ and is expressed by
\begin{equation}\label{05-09-2}
 p^*_{\alpha,\beta} = \frac{N}{N-2\alpha+1+\beta}.
\end{equation}
We observe that $p^*_{\alpha,0}=p^*_\alpha$ and $p^*_{\alpha,\beta}<p^*_{\alpha}$ when $\beta>0$.
In the source case, the assumptions on $g$ are of a different nature from in the absorption case, namely\\
\begin{itemize}
\item[$(G)\ $]
\begin{enumerate}\item[]
$g:\R_+\mapsto\R_+$
is a continuous  function  which satisfies
\begin{equation}\label{06-08-2}
g(s)\le c_1s^p+\sigma_0,\quad \forall s\ge0,
\end{equation}
for some $p\in(0,p^*_{\alpha,\beta})$, where $c_1>0$ and $\sigma_0>0$.
\end{enumerate}
\end{itemize}

Our main result concerning the source case is the following.

\begin{teo}\label{teo 10}
Assume that $\epsilon=-1$, $\nu\in\mathfrak{M}(\Omega,\rho^{\beta})$ with $\beta\in[0,2\alpha-1)$ is nonnegative,
 $g$ satisfies $(G)$ and\smallskip

\noindent (i) $p\in(0,1)$,
or\smallskip

\noindent (ii)
$p=1$ and $c_1$ is small enough,  or\smallskip

\noindent (iii)
$p\in(1,p^*_{\alpha,\beta})$, $\sigma_0$ and $\norm{\nu}_{\mathfrak{M}(\Omega,\rho^\beta)}$ are small enough.\smallskip

\noindent Then problem (\ref{eq1.1}) admits a  weak nonnegative solution $u_\nu$ which satisfies
\begin{equation}\label{1.5'}
 u_\nu\ge \mathbb{G}_\alpha[\nu].
\end{equation}
\end{teo}

We note that  Bidaut-V\'eron, Garc\'ia-Huidobro and V\'eron in \cite{BGV} obtained the existence of a renormalized solution of
$$-\Delta_p u=|\nabla u|^q+\nu\quad {\rm in }\ \Omega,$$
when $\nu\in \mathfrak{M}^b(\Omega)$. We make use of some idea  in \cite{BGV} in the proof of Theorem \ref{teo 10} and extend some results in \cite{BGV} to elliptic equations involving $(-\Delta)^\alpha$ with $\alpha\in(1/2,1)$ and   $\nu\in\mathfrak{M}(\Omega,\rho^\beta)$ with $\beta\in[0,2\alpha-1)$.

In the last section, we assume that $\Omega$ contains $0$ and give pointwise estimates of the positive solutions
\begin{equation}\label{eqX.1}
 \arraycolsep=1pt
\begin{array}{lll}
 (-\Delta)^\alpha  u+|\nabla u|^p=\delta_0\quad & \rm{in}\quad\Omega,\\[2mm]
 \phantom{ (-\Delta)^\alpha  +|\nabla u|^p}
u=0\quad & \rm{in}\quad \Omega^c,
\end{array}
\end{equation}
when $0<p<p^*_\alpha$. Combining properties of the Riesz kernel  with a bootstrap argument, we prove that any weak solution of (\ref{eqX.1}) is regular outside $0$ and is actually a classical solution of
\begin{equation}\label{eqX.2}
 \arraycolsep=1pt
\begin{array}{lll}
 (-\Delta)^\alpha  u+|\nabla u|^p=0\quad & \rm{in}\quad\Omega\setminus\{0\},\\[2mm]
 \phantom{ (-\Delta)^\alpha  +|\nabla u|^p}
u=0\quad & \rm{in}\quad \Omega^c.
\end{array}
\end{equation}
These pointwise estimates are quite easy to establish in the case $\alpha=1$, but much more delicate when the diffusion operator is non-local. We give sharp asymptotics of the behaviour of $u$ near $0$ and prove that the solution of (\ref{eqX.1}) is unique in the class of positive solutions. \smallskip

The paper is organized as follows. In Section 2, we study the Green operator and prove the key estimate
$$
\|\nabla \mathbb{G}_\alpha[\nu]\|_{M^{p^*_\alpha}(\Omega,\rho^\alpha  dx)}\le
c_2\|\nu\|_{\mathfrak{M}(\Omega,\rho^\beta)}
$$ Section  3 is devoted to prove Theorem \ref{teo 1} and Theorem \ref{teo 10}.
In Section 4, we consider the case where $\epsilon=1$ in  (\ref{eq1.1}) and $\nu$ is a Dirac mass. We obtain precise asymptotic estimate and derive uniqueness.\smallskip

\noindent{\it Aknowledgements.}  The authors are grateful to Marie-Fran\c{c}oise Bidaut-V\'eron for useful discussions in the preparation of this work.

\setcounter{equation}{0}
\section{Preliminaries}

\subsection{Marcinkiewicz type estimates}

In this subsection,
we recall some definitions and  properties of Marcinkiewicz
spaces.

\begin{definition}
Let $\Theta\subset \R^N$ be a domain and $\mu$ be a positive
Borel measure in $\Theta$. For $\kappa>1$,
$\kappa'=\kappa/(\kappa-1)$ and $u\in L^1_{loc}(\Theta,d\mu)$, we
set
\begin{equation}\label{mod M}
\|u\|_{M^\kappa(\Theta,d\mu)}=\inf\left\{c\in[0,\infty]:\int_E|u|d\mu\le
c\left(\int_Ed\mu\right)^{\frac1{\kappa'}},\ \forall E\subset \Theta,\,E\
\rm{Borel}\right\}
\end{equation}
and
\begin{equation}\label{spa M}
M^\kappa(\Theta,d\mu)=\{u\in
L_{loc}^1(\Theta,d\mu):\|u\|_{M^\kappa(\Theta,d\mu)}<\infty\}.
\end{equation}
\end{definition}

$M^\kappa(\Theta,d\mu)$ is called the Marcinkiewicz space of
exponent $\kappa$, or weak $L^\kappa$-space and
$\|.\|_{M^\kappa(\Theta,d\mu)}$ is a quasi-norm.

\begin{proposition}\label{pr 1} \cite{BBC,CC}
Assume that $1\le q< \kappa<\infty$ and $u\in L^1_{loc}(\Theta,d\mu)$.
Then there exists  $c_3>0$ dependent of $q,\kappa$ such that
$$\int_E |u|^q d\mu\le c_3\|u\|_{M^\kappa(\Theta,d\mu)}\left(\int_E d\mu\right)^{1-q/\kappa},$$
for any Borel set $E$ of $\Theta$.
\end{proposition}

The next estimate is the key-stone in the proof of Theorem \ref{teo 1}.

\begin{proposition}\label{general}
Let  $\Omega\subset \R^N\ (N\ge2)$ be a  bounded $C^2$ domain
and $\nu\in\mathfrak{M}(\Omega,\rho^\beta)$ with $\beta\in[0,2\alpha-1]$. Then there exists $c_2>0$ such that
\begin{equation}\label{annex 0}
\|\nabla\mathbb{G}_\alpha[|\nu|]\|_{M^{p^*_\alpha}(\Omega,\rho^\alpha dx)}\le c_2\|\nu\|_{\mathfrak{M}(\Omega,\rho^\beta)},
\end{equation}
where  $\displaystyle\nabla\mathbb{G}_\alpha[|\nu|](x)=\int_{\Omega} \nabla_xG_\alpha(x,y)d
|\nu(y)|$ and $p^*_\alpha$ is given by (\ref{annex 00}).

\end{proposition}
{\it Proof}.  For $\lambda>0$ and $y\in \Omega$, we set
$$\omega_\lambda(y)=\left\{x\in\Omega\setminus\{y\}:|\nabla_x G_\alpha(x,y)|\rho^\alpha(x)>\lambda \right\},\;
 m_\lambda(y)=\int_{\omega_\lambda(y)}dx.$$
From \cite{CS}, there exists $c_4>0$ such that for any $(x,y)\in
\Omega\times\Omega$ with $x\neq y$,
\begin{equation}\label{annex 01}
G_\alpha(x,y)\le c_4
\min\left\{\frac1{|x-y|^{N-2\alpha}},\frac{\rho^\alpha(x)}{|x-y|^{N-\alpha}},\frac{\rho^\alpha(y)}{|x-y|^{N-\alpha}}\right\},
\end{equation}
$$
G_\alpha(x,y)\le c_4 \frac{\rho^\alpha(y)}{\rho^\alpha(x)|x-y|^{N-2\alpha}},
$$
and by Corollary 3.3 in \cite{BKN}, we have
\begin{equation}\label{cap7 annex 01}
 |\nabla_x G_\alpha(x,y)|\le N G_\alpha(x,y)\max\left\{\frac1{|x-y|},\frac1{\rho(x)}\right\}.
\end{equation}
This implies that for any $\tau\in[0,1]$
\begin{eqnarray*}
  G_\alpha(x,y)\le c_4(\frac{\rho^\alpha(y)}{|x-y|^{N-\alpha}})^\tau(\frac{\rho^\alpha(x)}{|x-y|^{N-\alpha}})^{1-\tau}
   = c_4 \frac{\rho^{\alpha\tau}(y)\rho^{\alpha(1-\tau)}(x)}{|x-y|^{N-\alpha}},
\end{eqnarray*}
and then
\begin{equation}\label{annex 0001}
 |\nabla_x G_\alpha(x,y)|\le  c_5\max\left\{\frac{\rho^{\alpha}(y)}{\rho^\alpha(x)|x-y|^{N-2\alpha+1}},
\frac{\rho^{\alpha\tau}(y)\rho^{\alpha(1-\tau)-1}(x)}{|x-y|^{N-\alpha}}\right\}.
\end{equation}
Letting $\tau=\frac{2\alpha-1}{\alpha}\frac{N-\alpha}{N-2\alpha+1}\in(0,1)$, we derive
\begin{eqnarray*}
|\nabla_x G_\alpha(x,y)|\rho^{\alpha}(x)\le  c_5\max\left\{\frac{\rho^{2\alpha-1}(y)\rho_\Omega^{1-\alpha}}{|x-y|^{N-2\alpha+1}}, \frac{\rho^{\frac{(2\alpha-1)(N-\alpha)}{N-2\alpha+1}}(y)\rho_\Omega^{\frac{(2\alpha-1)(1-\alpha)}{N-2\alpha+1}}}{|x-y|^{N-\alpha}}\right\}.
\end{eqnarray*}
where $\rho_{\Omega}=\sup_{z\in\Omega}\rho(z)$.
There exists some $c_6>0$ such that
$$\omega_\lambda(y)\subset \left\{x\in\Omega: |x-y|\le c_6\rho^{\frac{2\alpha-1}{N-2\alpha+1} }(y)\max\{ \lambda^{-\frac1{N-2\alpha+1}}, \lambda^{-\frac1{N-\alpha}}\}\right\}.$$
By  $N-2\alpha+1>N-\alpha$, we deduce that for any $\lambda>1$, there holds
\begin{equation}\label{gra 0003}
\omega_\lambda(y)\subset \{x\in\Omega: |x-y|\le c_6 \rho^{\frac{2\alpha-1}{N-2\alpha+1}}(y) \lambda^{-\frac1{N-2\alpha+1}}\}.
\end{equation}
As a consequence,
$$ m_\lambda(y)\le c_7\rho^{(2\alpha-1)p_\alpha^*}(y)\lambda^{-p_\alpha^*},$$
where $c_7>0$ independent of $y$ and $\lambda$.

Let $E\subset\Omega$  be a Borel set and $\lambda>1$, then
\begin{eqnarray*}
\int_E|\nabla_x G_\alpha(x,y)|\rho^\alpha(x) dx\le
\int_{\omega_\lambda(y)}|\nabla_x G_\alpha(x,y)|\rho^\alpha(x)dx+\lambda\int_Edx.
\end{eqnarray*}
Noting that
\begin{eqnarray*}
\int_{\omega_\lambda(y)}|\nabla_x G_\alpha(x,y)|\rho^\alpha(x)dx&=&-\int_{\lambda}^\infty s
dm_s(y)
\\&=&\lambda m_\lambda(y)+ \int_{\lambda}^\infty m_s(y)ds
\\&\le& c_8\rho^{(2\alpha-1)p_\alpha^*}(y)\lambda^{1-p_\alpha^*},
\end{eqnarray*}
for some $c_8>0$, we derive
\begin{eqnarray*}
\int_E |\nabla_x G_\alpha(x,y)|\rho^\alpha(x)dx\le c_8\rho^{(2\alpha-1)p^*_\alpha}(y)\lambda^{1-p_\alpha^*}+\lambda \int_Edx.
\end{eqnarray*}
Choosing $\lambda=\rho^{2\alpha-1}(y)(\int_Edx)^{-\frac1{p_\alpha^*}}$ yields
\begin{eqnarray*}
\int_E |\nabla_x G_\alpha(x,y)|\rho^\alpha(x)dx\le (c_8+1)\rho^{2\alpha-1}(y)(\int_E dx)^{\frac{p_\alpha^*-1}{p_\alpha^*}},\quad \forall y\in \Omega.
\end{eqnarray*}
Therefore,
\begin{equation}\label{none}\begin{array}{lll}\displaystyle
\displaystyle \int_E|\nabla\mathbb{G}_\alpha[|\nu|](x)|\rho^\alpha(x)dx=\int_\Omega\int_E
|\nabla_x G_\alpha(x,y)|\rho^\alpha(x) dx d|\nu(y)|
\\[4mm]\phantom{\int_E|\nabla\mathbb{G}_\alpha[|\nu|](x)|}
\displaystyle\le \int_\Omega\rho^{2\alpha-1}(y)\left(\rho^{1-2\alpha}(y)\int_E
|\nabla_x G_\alpha(x,y)|\rho^\alpha(x)dx \right)d|\nu(y)|
\\[4mm]\phantom{\int_E|\nabla\mathbb{G}_\alpha[|\nu|](x)|}
\displaystyle\le (c_8+1)\int_\Omega \rho^\beta(y)\rho^{2\alpha-1-\beta}(y)d|\nu(y)|\left(\int_E dx\right)^{\frac{p^*_\alpha-1}{p^*_\alpha}}
\\[4mm]\phantom{\int_E|\nabla\mathbb{G}_\alpha[|\nu|](x)|}
\displaystyle\le (c_8+1)\rho_\Omega^{2\alpha-1-\beta}\|\nu\|_{\mathfrak{M}(\Omega,\rho^\beta)} \left(\int_E
dx\right)^{\frac{p^*_\alpha-1}{p^*_\alpha}}.
\end{array}
\end{equation}
As a consequence,
\begin{eqnarray*}
\|\nabla\mathbb{G}_\alpha[|\nu|]\|_{M^{p^*_\alpha}(\Omega,\rho^\alpha  dx)}\le
c_2\|\nu\|_{\mathfrak{M}(\Omega,\rho^\beta)},
\end{eqnarray*}
which ends the proof.\hfill$\Box$

\begin{proposition}\label{pr 4}\cite{CV2} Assume that  $\nu\in L^1(\Omega,\rho^\beta dx)$ with $0\leq \beta\leq\alpha$. Then for $p\in (1,\frac{N}{N-2\alpha+\beta})$, there exists $c_9>0$ such that for any $\nu\in L^{1}(\Omega,\rho^{\beta}dx)$
  \begin{equation}\label{power1'}
  \norm{\mathbb G_\alpha[\nu]}_{W^{2\alpha-\gamma,p}(\Omega)}\leq c_9\norm \nu_{L^{1}(\Omega,\rho^{\beta}dx)},
\end{equation}
where $p'=\frac p{p-1}$, $\gamma=\beta+\frac{N}{p'}$ if $\beta>0$ and
$\gamma>\frac{N}{p'}$ if $\beta=0$.
 \end{proposition}

\begin{proposition}\label{pr5} If $0\leq \beta<2\alpha-1$, then the mapping $\nu\mapsto |\nabla\mathbb G_\alpha[\nu]|$ is compact from  $L^{1}(\Omega,\rho^\beta dx)$ into $L^{q}(\Omega)$ for any $q\in [1,p^*_{\alpha,\beta})$ and there exists  $c_{10}>0$ such that
\begin{equation}\label{05-09-3}
\left (\int_\Omega |\nabla\mathbb G_\alpha[\nu](x)|^q dx\right)^{\frac1q}\le c_{10}\int_\Omega |\nu(x)| \rho^\beta(x) dx,
\end{equation}
where $p_{\alpha,\beta}^*$ is given by (\ref{05-09-2}).
 \end{proposition}
  {\it Proof.} For $\nu\in L^1(\Omega,\rho^\beta dx)$ with $0\leq \beta<2\alpha-1<\alpha$ , we  obtain from Proposition \ref{pr 4} that
  $$\mathbb G_\alpha[\nu]\in W^{2\alpha-\gamma,p}(\Omega),$$
  where $p\in(1,p_{\alpha,\beta}^*)$ and $2\alpha-\gamma>1$.  Therefore, $|\nabla\mathbb G_\alpha[\nu]|\in W^{2\alpha-\gamma-1,p}(\Omega)$ and
   \begin{equation}\label{06-09-1010}
     \|\nabla\mathbb G_\alpha[\nu]\|_{W^{2\alpha-\gamma-1,p}(\Omega)}\leq c_9\norm \nu_{L^{1}(\Omega,\rho^{\beta}dx)}.
   \end{equation}
 By \cite[Corollary 7.2]{NPV}, the embedding of $W^{2\alpha-\gamma-1,p}(\Omega)$ into $L^{q}(\Omega)$ is compact for
 $q\in[1,\frac{Np}{N-(2\alpha-\gamma-1)p})$.  When $\beta>0$,
\begin{eqnarray*}
\frac{Np}{N-(2\alpha-\gamma-1)p} &=& \frac{Np}{N-(2\alpha-\beta-N\frac{p-1}{p}-1)p} \\
   &=& \frac{N}{N-2\alpha+1+\beta}=p^*_{\alpha,\beta}.
\end{eqnarray*}
When $\beta=0$,
\begin{eqnarray*}
\lim_{\gamma\to(\frac{N}{p'})^+}\frac{Np}{N-(2\alpha-\gamma-1)p} &=& \frac{Np}{N-(2\alpha-N\frac{p-1}{p}-1)p} \\
   &=& \frac{N}{N-2\alpha+1}=p^*_{\alpha,0}.
\end{eqnarray*}
Then the mapping $\nu\mapsto |\nabla\mathbb G_\alpha[\nu]|$ is compact from  $L^{1}(\Omega,\rho^\beta dx)$ into $L^{q}(\Omega)$ for any $q\in [1,p^*_{\alpha,\beta})$. Inequality (\ref{05-09-3}) follows by  (\ref{06-09-1010}) and the continuity of the embedding  of $W^{2\alpha-\gamma-1,p}(\Omega)$ into $L^{q}(\Omega)$. \hfill$\Box$\medskip

\noindent{\it Remark.}
If  $\nu\in L^1(\Omega,\rho^\beta dx)$ with $0\leq \beta<2\alpha-1$ and $u$ is the  solution of
$$\arraycolsep=1pt
\begin{array}{lll}
 (-\Delta)^\alpha  u=\nu\quad & \rm{in}\quad\Omega,\\[2mm]
 \phantom{   (-\Delta)^\alpha  }
u=0\quad & \rm{in}\quad \Omega^c,
\end{array}
$$
then for any $q\in[1,p^*_{\alpha,\beta})$,
$$\left(\int_\Omega|\nabla u|^qdx\right)^{\frac1q}\le c_{10}\int_\Omega|\nu(x)|\rho^\beta(x) dx.$$

\subsection{Classical solutions }

In this subsection we consider the question of existence of classical solutions to problem
\begin{equation}\label{5.4}
 \arraycolsep=1pt
\begin{array}{lll}
 (-\Delta)^\alpha  u+h(|\nabla u|)=f\quad & \rm{in}\quad\Omega,\\[2mm]
 \phantom{   (-\Delta)^\alpha  u+g(|\nabla u|}
u=0\quad & \rm{in}\quad \Omega^c.
\end{array}
\end{equation}
%

\begin{teo}\label{pr 2} Assume $h\in C^\theta(\R_+)\cap L^{\infty}(\R_+)$ for some $\theta\in (0,1]$
 and  $f\in C^\theta(\bar\Omega)$. Then problem (\ref{5.4})
 admits a unique classical solution $u$. Moreover,

\noindent $(i)$  if $f-h(0)\ge0$ in $\Omega$, then $u\ge0$; \smallskip

\noindent $(ii)$ the mappings $h\mapsto u$ and $f\mapsto u$ are respectively nonincreasing and nondecreasing.

\end{teo}

\noindent{\it Proof.} We divide the proof into several steps.\\
 {\it Step 1. Existence}. We define the operator $T$ by
 $$Tu=\mathbb{G}_\alpha\left[f-h(|\nabla u|)\right],\quad \forall u\in W^{1,1}_0(\Omega).$$
Using (\ref{annex 0001}) with $\tau=0$ yields
\begin{eqnarray}
\norm {Tu}_{W^{1,1}(\Omega)}&\leq&\nonumber
\norm{\mathbb{G}_\alpha [f]}_{W^{1,1}(\Omega)}+\norm{\mathbb{G}_\alpha [h(|\nabla u|)]}_{W^{1,1}(\Omega)}
\\&\le&\left(\norm{f}_{L^\infty(\Omega)}+\norm{h(|\nabla u|)}_{L^\infty(\Omega)}\right)\norm{\int_\Omega G_\alpha(\cdot,y)dy }_{W^{1,1}(\Omega)}\nonumber
\\&=& c_{11}\left(\norm{f}_{L^\infty(\Omega)}+\norm{h}_{L^\infty(\R_+)}\right),\label{2.2}
\end{eqnarray}
where  $c_{11}=\norm{\int_\Omega G_\alpha(\cdot,y)dy }_{W^{1,1}(\Omega)}$. Thus $T$ maps $W_0^{1,1}(\Omega)$ into itself. Clearly, if $u_n\to u$ in $W_0^{1,1}(\Omega)$ as $n\to\infty$, then $h(|\nabla u_n|)\to h(|\nabla u|)$ in $L^1(\Omega)$, thus $T$ is continuous.
We claim that $T$ is a compact operator.
In fact, for $u\in W^{1,1}_0(\Omega)$, we see that
$f-h(|\nabla u|)\in L^1(\Omega)$ and then, by Proposition \ref{pr 4},
it implies that $Tu\in W_0^{2\alpha-\gamma,p}(\Omega)$
where $\gamma\in(\frac{N(p-1)}{p},2\alpha-1)$ and $2\alpha-1>\frac{N(p-1)}{p}>0$ for $p\in (1,\frac N{N-2\alpha+1})$.
Since  the embedding $W_0^{2\alpha-\gamma,p}(\Omega)\hookrightarrow W^{1,1}_0(\Omega)$ is  compact, $T$ is a compact operator.

Let $\mathcal{O}=\{u\in W^{1,1}_0(\Omega): \norm {u}_{W^{1,1}(\Omega)}\le c_{10}(\norm{f}_{L^\infty(\Omega)}+\norm{h}_{L^\infty(\R_+)}) \}$, which is a closed and convex
set of $W^{1,1}_0(\Omega)$.  Combining with (\ref{2.2}), there holds
$$T(\mathcal{O})\subset \mathcal{O}.$$
 It follows by Schauder's fixed point theorem that there exists some $u\in W_0^{1,1}(\Omega)$ such that
$Tu=u$.

Next we show that $u$ is a classical solution of  (\ref{5.4}). Let  open set $O$ satisfy $ O\subset \bar O\subset \Omega$.
By Proposition 2.3 in \cite{RS}, for any $\sigma\in(0,2\alpha)$, there exists $c_{12}>0$ such that
$$\norm{u}_{C^{\sigma}(O)}\le c_{12}\{\norm{h(|\nabla u|)}_{L^\infty(\Omega)}+\norm{f}_{L^{\infty}(\Omega)}\},$$
and by choosing $\sigma=\frac{2\alpha+1}2\in(1,2\alpha)$, then
$$\norm{|\nabla u|}_{C^{\sigma-1}(O)}\le c_{12}\{\norm{h(|\nabla u|)}_{L^\infty(\Omega)}+\norm{f}_{L^{\infty}(\Omega)}\},$$
and then applied  \cite[Corollary 2.4]{RS}, $u$ is $C^{2\alpha+\epsilon_0}$ locally in $\Omega$ for some $\epsilon_0>0$.
Then $u$ is a classical solution of (\ref{5.4}).
Moreover, from \cite{CV2}, we have
\begin{equation}\label{5.6}
\int_\Omega [u(-\Delta)^\alpha\xi+h(|\nabla u|)\xi]dx=\int_\Omega\xi
fdx,\quad \forall\xi\in \mathbb{X}_{\alpha}.
\end{equation}

\noindent{\it Step 2. Proof of $(i)$.} If $u$ is not nonnegative,  then there exists $x_0\in\Omega$ such that
$$u(x_0)=\min_{x\in\Omega} u(x)< 0,$$
 then $\nabla u(x_0)=0$ and $(-\Delta)^\alpha u(x_0)<0$.
Since $u$ is the classical solution of  (\ref{5.4}), $(-\Delta)^\alpha u(x_0)=f(x_0)-h(0)\ge0$,
which is a contradiction.

\smallskip

\noindent{\it Step 3. Proof of $(ii)$.} We just give the proof of the first argument,
the proof of the second being similar.
Let $h_1$ and $h_2$ satisfy our hypotheses for $h$ and
$h_1\le h_2$. Denote $u_1$ and $u_2$ the solutions of (\ref{5.4}) with $h$ replaced by $h_1$ and $h_2$ respectively.
If there exists
 $x_0\in \Omega$ such that
 $$(u_1-u_2)(x_0)=\min_{x\in\Omega}\{(u_1-u_2)(x) \}<0. $$
Then
$$(-\Delta)^\alpha(u_1-u_2)(x_0)<0,\quad \nabla u_1(x_0)=\nabla u_2(x_0).$$
This implies
\begin{eqnarray}\label{5.7}
(-\Delta)^\alpha (u_1-u_2)(x_0)+h_1(|\nabla u_1(x_0)|)-h_2(|\nabla u_2(x_0)|)<0.
\end{eqnarray}
However,
\begin{eqnarray*}
(-\Delta)^\alpha (u_1-u_2)(x_0)+h_1(|\nabla u_1(x_0)|)-h_2(|\nabla u_2(x_0)|)=f(x_0)-f(x_0)=0,
\end{eqnarray*}
contradiction. Then $u_1\ge u_2$.
\smallskip

\noindent Uniqueness follows from Step 3.\hfill$\Box$

\setcounter{equation}{0}
\section{Proof of Theorems \ref{teo 1} and \ref{teo 10}}

\subsection{The absorption case}
In this subsection, we prove the existence of a weak solution to (\ref{eq1.1}) when
$\epsilon=1$. To this end, we give below an auxiliary lemma.
\begin{lemma}\label{lm 08-09}
Assume that $g:\R_+\mapsto\R_+$ is  continuous  and (\ref{1.4}) holds with $p^*_{\alpha}$. Then there is a sequence
real positive numbers $\{T_n\}$ such that
$$\lim_{n\to\infty}T_n=\infty\quad{\rm and}\quad \lim_{n\to\infty}g(T_n)T_n^{-p^*_{\alpha}}=0.$$

\end{lemma}
{\it Proof.} Let $\{s_n\}$ be a sequence of real positive numbers converging to $\infty$. We observe
\begin{eqnarray*}
\int_{s_n}^{2s_n}g(t)t^{-1-p^*_{\alpha}}dt&\ge&
\min_{t\in[s_n,2s_n]}g(t)(2s_n)^{-1-p^*_{\alpha}}\int_{s_n}^{2s_n}dt
\\&=&2^{-1-p^*_{\alpha}}s_n^{-p^*_{\alpha}}\min_{t\in[s_n,2s_n]}g(t)
\end{eqnarray*}
and by (\ref{1.4}),
\begin{eqnarray*}
\lim_{n\to\infty}\int_{s_n}^{2s_n}g(t)t^{-1-p^*_{\alpha}}dt=0.
\end{eqnarray*}
Then we choose $T_n\in[s_n,2s_n]$ such that $g(T_n)=\min_{t\in[s_n,2s_n]}g(t)$ and then
the claim follows.
\hfill$\Box$
\medskip

\noindent{\it Proof of Theorem \ref{teo 1}.}
Let $\beta\in[0,2\alpha-1)$, we define the space
$$C_{\beta}(\bar \Omega)=\{\zeta\in C(\bar \Omega):\rho^{-\beta}\zeta\in C(\bar \Omega)\}$$
endowed with the norm
$$\norm{\zeta}_{C_{\beta}(\bar\Omega)}=\|\rho^{-\beta}\zeta\|_{C(\bar\Omega)}. $$
Let  $\{\nu_n\}\subset C^1(\bar \Omega)$ be a sequence of nonnegative functions such that
$\nu_{n }\to\nu $ in  sense of duality with $C_{\beta}(\bar
\Omega)$, that is,
\begin{equation}\label{06-08}
  \lim_{n\to\infty}\int_{\bar \Omega}\zeta \nu_{n }dx=\int_{\bar \Omega}\zeta d\nu,\qquad\forall \zeta\in C_{\beta}(\bar \Omega).
\end{equation}
By the
Banach-Steinhaus Theorem, $\norm{\nu_{n}}_{\mathfrak M
(\Omega,\rho^\beta)}$ is bounded independently of $n$. We consider a sequence $\{g_n\}$ of $C^1$ nonnegative  functions defined on $\R_+$
such that $g_n(0)=0$ and
\begin{equation}\label{06-08-1}
  g_n\le g_{n+1}\le g,\quad \sup_{s\in\R_+}g_n(s)=n\quad{\rm and}\quad \lim_{n\to\infty}\norm{g_n-g}_{L^\infty_{loc}(\R_+)}=0.
\end{equation}
By Theorem \ref{pr 2}, there exists a unique \emph{nonnegative} solution $u_n$ of (\ref{eq1.1}) with data $\nu_n$ and $g_n$ instead of $\nu$ and $g$, and there holds
\begin{equation}\label{L7}
\displaystyle\int_\Omega\left( u_n +g_n(|\nabla u_n|)\eta_1\right)dx
=\displaystyle \int_\Omega \nu_n\eta_1dx\leq C\norm{\nu}_{\mathfrak M (\Omega,\rho^\beta)},
\end{equation}
where $\eta_1=\mathbb{G}_\alpha[1]$.
Therefore, $\norm{g_n(|\nabla u_n|)}_{\mathfrak M (\Omega,\rho^\alpha)}$ is
bounded independently of $n$. For $\varepsilon>0$ and
$\xi_{\varepsilon}=(\eta_1+\varepsilon)^{\frac{\beta}{\alpha}}-\varepsilon^{\frac{\beta}{\alpha}}\in \mathbb{X}_{\alpha}$
which is concave in the interval $[0,\eta_1(\bar\omega)]$,
 where $\eta_1(\bar\omega)=\max_{x\in\Omega}\eta_1(x)$. By  \cite[Lemma 2.3 (ii)]{CV2}, we see that
$$\displaystyle\begin{array}{lll}\displaystyle
(-\Delta)^{\alpha}\xi_{\varepsilon}=\frac{\beta}{\alpha}(\eta_1+\varepsilon)^{\frac1{\alpha}}
(-\Delta)^{\alpha}\eta_1-\frac{\beta(\beta-\alpha)}{\alpha^2}(\eta_1+\varepsilon)^{\frac{\beta-2\alpha}{\alpha}}
\int_{\Omega}\frac{(\eta_1(y)-\eta_1(x))^2}{|y-x|^{N+2\alpha}}dy\\[4mm]
\phantom{(-\Delta)^{\alpha}\xi_{\varepsilon}}\displaystyle \geq
\frac{\beta}{\alpha}(\eta_1+\varepsilon)^{\frac{\beta-\alpha}{\alpha}},
\end{array}$$
and $\xi_\varepsilon\in\mathbb{X}_{\alpha}$. Since
$$\int_{\Omega}\left(u_n(-\Delta)^{\alpha}\xi_{\varepsilon}+g_n(|\nabla u_n|)\xi_{\varepsilon}\right)dx=
\int_{\Omega}\xi_{\varepsilon} \nu_n dx,
$$
we obtain
$$\int_{\Omega}\left(\frac{\beta}{\alpha}u_n(\eta_1+\varepsilon)^{\frac{\beta-\alpha}{\alpha}}+g_n(|\nabla u_n|)\xi_{\varepsilon}\right)dx\leq
\int_{\Omega}\xi_{\varepsilon} \nu_n dx.
$$
If we let $\varepsilon\to 0$, it yields
$$\int_{\Omega}\left(\frac{\beta}{\alpha}u_n\eta_1^{\frac{\beta-\alpha}{\alpha}}+g_n(|\nabla u_n|)\eta_1^{\frac{\beta}{\alpha}}\right)dx\leq
\int_{\Omega}\eta_1^{\frac{\beta}{\alpha}} \nu_n dx.
$$
Using \cite[Lemma 2.3]{CV2}, we derive the estimate
\begin{equation}\label{L9}
\int_{\Omega}\left(u_n\rho^{\beta-\alpha}+g_n(|\nabla u_n|)\rho^{\beta}\right)
dx\leq c_{13}\norm {\nu_n}_{\mathfrak M (\Omega,\rho^{\beta})}\leq c_{14}\norm {\nu}_{\mathfrak M (\Omega,\rho^{\beta})}.
\end{equation}
Thus $\{g_n(|\nabla u_n|)\}$ is uniformly bounded in $L^1(\Omega,\rho^{\beta}dx).$
Since $u_n=\mathbb G[\nu_n-g_n(|\nabla u_n|)]$, there holds
\begin{eqnarray*}
\norm {|\nabla u_n|}_{M^{p^*_\alpha}(\Omega,\rho^\alpha dx)}&\leq&
\norm{\nu_{n}}_{\mathfrak M (\Omega,\rho^{\beta})}+\norm{g_n(|\nabla u_n|)}_{\mathfrak M (\Omega,\rho^{\beta})}
\\&\le& c_{15}\norm{\nu}_{\mathfrak M (\Omega,\rho^{\beta})}.
\end{eqnarray*}
Since $\nu_n-g_n(|\nabla u_n|)$ is uniformly bounded in $L^1(\Omega,\rho^{\beta}dx)$, we use Proposition \ref{pr5} to obtain that the sequences $\{ u_n\}$, $\{|\nabla u_n|\}$ are relatively compact in
$L^q(\Omega)$ for $ q\in[1,\frac{N}{N-2\alpha+\beta})$ and $q\in[1,p^*_{\alpha,\beta})$, respectively. Thus,
there exist a sub-sequence $\{u_{n_k}\}$ and some $u\in L^q(\Omega)$ with  $q\in[1,\frac{N}{N-2\alpha+\beta})$
  such that \smallskip

\noindent $(i)$ $ u_{n_k}\to u$ a.e. in $\Omega$ and in $L^q(\Omega)$ with $ q\in[1,\frac{N}{N-2\alpha+\beta})$;

\noindent $(ii)$  $|\nabla  u_{n_k}|\to |\nabla u|$ a.e. in $\Omega$ and in $L^q(\Omega)$ with $q\in[1,p^*_{\alpha,\beta})$.\smallskip

\noindent Therefore,
$g_{n_k}(|\nabla u_{n_k}|)\to g(|\nabla u|)$ a.e. in $\Omega$.   For $\lambda>0$, we denote
$$
S_\lambda=\{x\in\Omega:|\nabla u_{n_k}(x)|>\lambda\}\quad{\rm and }\quad \omega(\lambda)=\int_{S_\lambda}\rho^{\alpha}(x)dx.$$
Then for any Borel
set $E\subset\Omega$, we have that
$$\begin{array} {ll}\displaystyle
\int_{E}g_{n_k}(|\nabla u_{n_k}|)|\rho^{\alpha}(x)dx\le \int_{E}g(|\nabla u_{n_k}|)|\rho^{\alpha}(x)dx
\\[4mm]\displaystyle\phantom{-----}
 =\int_{E\cap S^c_{\lambda}}g(|\nabla u_{n_k}|)\rho^{\alpha}(x)dx+\int_{E\cap S_{\lambda}}g(|\nabla u_{n_k}|)\rho^{\alpha}(x)dx
 \\[4mm]\displaystyle\phantom{---------}\displaystyle
\le \tilde g(\lambda)\int_E\rho^{\alpha}(x)dx+\int_{S_{\lambda}}g(|\nabla u_{n_k}|)\rho^{\alpha}(x)dx
 \\[4mm]\displaystyle\phantom{-------------}\displaystyle
 \le \tilde g(\lambda)\int_E\rho^{\alpha}(x)dx-\int_{\lambda}^\infty g(s)d\omega(s),
\end{array}$$
where $\tilde g(s)=\max_{t\in[0,s]}\{g(t)\}$.
But
$$\int_{\lambda}^\infty  g(s)d\omega(s)=\lim_{n\to\infty}\int_{\lambda}^{T_n}  g(s)d\omega(s).
$$
where $\{T_n\}$ is given by Lemma \ref{lm 08-09}.
Since $|\nabla u_{n_k}|\in M^{p^*_\alpha}(\Omega,\rho^\alpha dx)$,
$\omega(s)\leq c_{16}s^{-p^*_\alpha}$ and
$$\displaystyle\begin{array}{lll}
\displaystyle-\int_{\lambda}^{T_n}   g(s)d\omega(s) =-\left[
g(s)\omega(s)\!\!\!\!\!\!\!\!^{\phantom{\frac{X^X}{X}}}\right]_{s=\lambda}^{s={T_n}}+\int_{\lambda}^{T_n}
\omega(s)d  g(s)
\\[4mm]\phantom{\int_{\lambda}^{T_n}   g(s)d\omega(s)}\displaystyle
\leq  g(\lambda)\omega(\lambda)- g({T_n})\omega({T_n})+c_{16}\int_{\lambda}^{T_n} s^{-p^*_\alpha}d  g(s)
\\[4mm]\phantom{\int_{\lambda}^{T_n}   g(s)d\omega(s)}\displaystyle
\leq   g(\lambda)\omega(\lambda)-  g({T_n})\omega({T_n})+
c_{16}\left({T_n}^{-p^*_\alpha}g({T_n})-\lambda^{-p^*_\alpha} g(\lambda)\right)
\\[4mm]\phantom{-----------\int_{\lambda}^{T_n}   g(s)d\omega(s)}\displaystyle
+\frac{c_{16}}{p^*_\alpha+1}\int_{\lambda}^{T_n}
s^{-1-p^*_\alpha}  g(s)ds.
\end{array}$$
 By assumption $(\ref{1.4})$ and Lemma \ref{lm 08-09},  it follows
\begin{equation}\label{07-09-20}
\lim_{n\to\infty}T_n^{-p^*_\alpha}g(T_n)=0.
\end{equation}
Along with $g(\lambda)\omega(\lambda)\leq
c_{16}\lambda^{-p^*_\alpha}  g(\lambda)$,  we have
$$-\int_{\lambda}^\infty   g(s)d\omega(s)\leq \frac{c_{16}}{p^*_\alpha+1}\int_{\lambda}^\infty s^{-1-p^*_\alpha} g(s)ds.
$$
Notice that the above quantity on the right-hand side tends to $0$
when $\lambda\to\infty$. It implies that for any
$\epsilon>0$ there exists $\lambda>0$ such that
$$\frac{c_{16}}{p^*_\alpha+1}\int_{\lambda}^\infty s^{-1-p^*_\alpha} g(s)ds\leq \frac{\epsilon}{2},
$$
and $\delta>0$ such that
$$\int_E \rho^\alpha(x) dx\leq \delta\Longrightarrow   \tilde g(\lambda)\int_E dx\leq\frac{\epsilon}{2}.
$$
This proves that $\{g_{n_k}( |\nabla u_{n_k}|)\}$ is uniformly integrable in
$L^1(\Omega,\rho^\alpha dx)$. Then $g_{n_k}(|\nabla u_{n_k}|)\to g( |\nabla u|)$ in
$L^1(\Omega,\rho^\alpha dx)$ by Vitali convergence theorem. Letting
$n_k\to\infty$ in the identity
$$\int_{\Omega}\left(u_{n_k}(-\Delta)^{\alpha}\xi+ g_{n_k} (|\nabla u_{n_k}|)\xi\right) dx=\int_{\Omega}\nu_{n_k}\xi dx,\quad \forall\xi\in \mathbb X_{\alpha},
$$
 it infers that $u$ is a weak
solution of $(\ref{eq1.1})$. Since $u_{n_k}$ is nonnegative, so is $u$.

\smallskip

\noindent Estimate (\ref{1.5}) is a consequence of  positivity and
\begin{eqnarray*}
 u_{n_k}=\mathbb{G}_\alpha[\nu_{n_k}]-\mathbb{G}_\alpha[g_{n_k}(|\nabla u_{n_k}|)]\le \mathbb{G}_\alpha[\nu_{n_k}].
\end{eqnarray*}
Since $\lim_{n_k\to\infty} u_{n_k}=u$, (\ref{1.5}) follows.
 \hfill$\Box$

\subsection{The source case}

In this subsection we study the existence of solutions to
problem (\ref{eq1.1}) when $\epsilon=-1$.

\smallskip

\noindent{\it Proof of Theorem \ref{teo 10}.} Let $\{\nu_n\}$ be a sequence of $C^2$ nonnegative functions converging to $\nu$ in the sense of
(\ref{06-08}), $\{g_n\}$ an increasing sequence of $C^1$, nonnegative bounded functions defined on $\R_+$ satisfying (\ref{06-08-1}) and converging to $g$.
We set $p_0=\frac{p+p^*_{\alpha,\beta}}2\in(p,p^*_{\alpha,\beta})$, where $p^*_{\alpha,\beta}$ is given by
(\ref{05-09-2}) and $p<p^*_{\alpha,\beta}$ is the maximal growth rate of $g$ which satisfies  (\ref{06-08-2}), and
$$M(v)=\left(\int_{\Omega} |\nabla v|^{p_0} dx\right)^{\frac1{p_0}}.$$
We may assume that $\norm{\nu_n}_{L^1(\Omega,\rho^{\beta}dx)}\le 2\norm{\nu}_{\mathfrak M (\Omega,\rho^{\beta})}$ for all $n\ge1$.\\
{\it Step 1. We claim that for $n\geq 1$,
$$ \arraycolsep=1pt
\begin{array}{lll}
 (-\Delta)^\alpha u_{n}=g_{n}(|\nabla u_{n}|)+\nu_n\quad & \rm{in}\quad\Omega,\\[2mm]
 \phantom{   (-\Delta)^\alpha_n }
u_{n}=0\quad & \rm{in}\quad \Omega^c
\end{array}
$$
admits a solution $u_n$ such that
$$M(u_n)\le \bar\lambda,$$
where $\bar\lambda>0$ independent of $n$. }

To this end, we define the operators $\{T_n\}$ by
 $$T_nu=\mathbb{G}_\alpha\left[g_n(|\nabla u|)+\nu_n\right],\quad \forall u\in W^{1,p_0}_0(\Omega).$$
On the one hand, using (\ref{annex 0001}) with $\tau=0$ yields
\begin{eqnarray*}
\norm {T_nu}_{W^{1,1}(\Omega)}&\leq&\nonumber
\norm{\mathbb{G}_\alpha [\nu_n]}_{W^{1,1}(\Omega)}+\norm{\mathbb{G}_\alpha [g_n(|\nabla u|)]}_{W^{1,1}(\Omega)}
\\&\le& c_{11}\left(\norm{\nu_n}_{L^\infty(\Omega)}+\norm{g_n}_{L^\infty(\R_+)}\right),
\end{eqnarray*}
where  $c_{11}=\norm{\int_\Omega G_\alpha(\cdot,y)dy }_{W^{1,1}(\Omega)}$.
On the other hand,  by (\ref{06-08-2}) and Proposition \ref{pr5}, we have
\begin{eqnarray}
  \left(\int_{\Omega}|\nabla (T_nu)|^{p_0} dx\right)^{\frac1{p_0}}&\le& c_2\norm{g_n(|\nabla u|)+\nu_n}_{L^1 (\Omega,\rho^{\beta}dx)}\nonumber\\
   &\le & c_2 [\norm{g_n(|\nabla u|)}_{L^1(\Omega,\rho^{\beta}dx)}+2\norm{\nu}_{\mathfrak M (\Omega,\rho^{\beta})}]\label{06-08-10}\\
   &\le & c_{2}c_1\int_{\Omega}|\nabla u|^p\rho^\beta dx+c_{17}\sigma_0 +2c_2\norm{\nu}_{\mathfrak M (\Omega,\rho^{\beta})},\nonumber
\end{eqnarray}
where $c_{17}=c_2\int_\Omega\rho^\beta dx$. Then we use H\"{o}lder inequality to obtain that
\begin{eqnarray}\label{06-08-11}
   \left(\int_{\Omega}|\nabla u|^p\rho^\beta dx\right)^{\frac1p}\le (\int_\Omega \rho^{\frac{\beta p_0}{p_0-p}}dx)^{\frac 1p-\frac1{p_0}}   \left(\int_{\Omega}|\nabla u|^{p_0} dx\right)^{\frac1{p_0}},
\end{eqnarray}
where $\int_\Omega \rho^{\frac{\beta p_0}{p_0-p}}dx$ is bounded, since $\frac{\beta p_0}{p_0-p}\ge0$.
Along with (\ref{06-08-10}) and (\ref{06-08-11}), we derive
\begin{equation}\label{05-09-4}
  M(T_nu)\le c_{18} M(u)^{p}+c_{19}\norm{\nu}_{\mathfrak M (\Omega,\rho^{\beta})}+c_{17}\sigma_0,
\end{equation}
where $c_{18}=c_2c_1(\int_\Omega \rho^{\frac{\beta p_0}{p_0-p}}dx)^{\frac 1p-\frac1{p_0}}>0$ and $c_{19}>0$ independent of $n$. Therefore, if we assume that $M(u)\le \lambda$, inequality (\ref{05-09-4}) implies
\begin{equation}\label{05-09-5}
  M(T_nu)\le c_{18}\lambda^p+c_{19}\norm{\nu}_{\mathfrak M (\Omega,\rho^{\beta})}+c_{17}\sigma_0.
\end{equation}
Let $\bar\lambda>0$ be the largest root of the equation
\begin{equation}\label{07-08-1}
c_{18}\lambda^{p}+c_{19}\norm{\nu}_{\mathfrak M (\Omega,\rho^{\beta})}+c_{17}\sigma_0= \lambda,
\end{equation}
This root exists if one of the following condition holds:\smallskip

\noindent (i) $p\in (0,1)$, in which case (\ref{07-08-1}) admits only one root;\smallskip

\noindent (ii) $p=1$ and $c_{17}<1$, and again (\ref{07-08-1}) admits only one root;\smallskip

\noindent (iii) $p\in (1,p^*_\alpha)$ and there exists $\varepsilon_0>0$ such that
$\max\left\{\norm{\nu}_{\mathfrak M (\Omega,\rho^{\beta})},\sigma_0\right\}\leq \varepsilon_0$.
In that case (\ref{07-08-1}) admits usually two positive roots.\smallskip

\noindent If we suppose that one of the above conditions holds, the definition of $\bar\lambda>0$ implies that it is the largest $\lambda>0$ such that
\begin{equation}\label{07-08-1'}
c_{18}\lambda^{p}+c_{19}\norm{\nu}_{\mathfrak M (\Omega,\rho^{\beta})}+c_{17}\sigma_0\leq  \lambda,
\end{equation}
For  $M(u)\le \bar\lambda$, we obtain that
  $$M(T_nu)\le c_{18}\bar\lambda^p+c_{19}\norm{\nu}_{\mathfrak M (\Omega,\rho^{\beta})}+c_{17}\sigma_0= \bar\lambda.
  $$
By the assumptions of Theorem \ref{teo 10}, $\bar\lambda$ exists and it is larger than $M(u_n)$. Therefore,
\begin{equation}\label{04-09}
 \int_{\Omega} |\nabla (T_nu)|^{p_0} dx\le \bar\lambda^{p_0}.
\end{equation}

Thus $T_n$ maps $W_0^{1,p_0}(\Omega)$ into itself. Clearly, if $u_n\to u$ in $W_0^{1,p_0}(\Omega)$ as $n\to\infty$, then $g_n(|\nabla u_n|)\to g_n(|\nabla u|)$ in $L^1(\Omega)$, thus $T$ is continuous.
We claim that $T$ is a compact operator.
In fact, for $u\in W^{1,p_0}_0(\Omega)$, we see that
$\nu_n-g_n(|\nabla u|)\in L^1(\Omega)$ and then, by Proposition \ref{pr 4},
it implies that $T_nu\in W_0^{2\alpha-\gamma,p}(\Omega)$
where $\gamma\in(\frac{N(p-1)}{p},2\alpha-1)$ and $2\alpha-1>\frac{N(p-1)}{p}>0$ for $p\in (1,\frac N{N-2\alpha+1})$.
Since  the embedding $W_0^{2\alpha-\gamma,p}(\Omega)\hookrightarrow W^{1,p_0}_0(\Omega)$ is  compact, $T_n$ is a compact operator.

Let
$$
\displaystyle\begin{array}{lll}\displaystyle
\mathcal{G}=\{u\in W^{1,p_0}_0(\Omega): \norm {u}_{W^{1,1}(\Omega)}\le c_{11}(\norm{\nu_n}_{L^\infty(\Omega)}+\norm{g_n}_{L^\infty(\R_+)})
\\[2mm]\phantom{---------}
{\rm and}\quad   M(u)\le \bar\lambda \},
\end{array}
$$
  which is a closed and convex
set of $W^{1,p_0}_0(\Omega)$.  Combining with (\ref{2.2}), there holds
$$T_n(\mathcal{G})\subset \mathcal{G}.$$
 It follows by Schauder's fixed point theorem that there exists some $u_n\in W_0^{1,p_0}(\Omega)$ such that
$T_nu_n=u_n$ and $M(u_n)\le \bar\lambda,$
where $\bar\lambda>0$ independent of $n$.
By the same arguments as in Theorem \ref{pr 2}, $u_n$ belongs to $C^{2\alpha+\epsilon_0}$ locally in $\Omega$, and
\begin{equation}\label{05-09-1}
 \int_\Omega u_{n}(-\Delta)^\alpha\xi=\int_\Omega g_{n}(|\nabla u_{n}|)\xi dx+\int_\Omega\xi\nu_ndx,\quad \forall \xi\in\mathbb{X}_\alpha.
\end{equation}

\noindent{\it Step 2: Convergence. }By (\ref{04-09}) and (\ref{06-08-11}),  $g_n(|\nabla u_n|)$ is uniformly bounded in $L^1(\Omega,\rho^\beta dx)$.
By Proposition \ref{pr 4},  $\{u_n\}$ is bounded in $W_0^{2\alpha-\gamma,q}(\Omega)$ where $q\in(1,p^*_{\alpha,\beta})$ and $2\alpha-\gamma>1$. By Proposition \ref{pr5},  there exist a subsequence $\{u_{n_k}\}$ and $u$ such that
$u_{n_k}\to u$ a.e. in $\Omega$ and in $L^1(\Omega)$, and
$|\nabla u_{n_k}| \to |\nabla u|$ a.e. in $\Omega$ and in $L^q(\Omega)$ for any $q\in[1,p^*_{\alpha,\beta})$.
By assumption (G), $g_{n_k}(|\nabla u_{n_k}|)\to g( |\nabla u|)$ in $L^1(\Omega)$.
Letting $n_k\to \infty$ to have that
 $$\int_\Omega u(-\Delta)^\alpha\xi=\int_\Omega g(|\nabla u|)\xi dx+\int_\Omega\xi d\nu,\quad \forall \xi\in\mathbb{X}_\alpha, $$
thus $u$ is a weak solution of (\ref{eq1.1}) which is nonnegative as $\{u_n\}$ are nonnegative. Furthermore, (\ref{1.5'}) follows from the positivity of $g(|\nabla u_n])$.
\hfill$\Box$

\setcounter{equation}{0}
\section{The case of the Dirac mass}

In this section we assume that $\Omega$ is an open, bounded and $C^2$ domain containing $0$ and  $u$ a nonnegative weak solution of
\begin{equation}\label{eq4.1}
 \arraycolsep=1pt
\begin{array}{lll}
 (-\Delta)^\alpha  u+ |\nabla u|^p=\delta_0\quad & \rm{in}\quad\Omega,\\[2mm]
 \phantom{   (-\Delta)^\alpha +|\nabla u|^p}
u=0\quad & \rm{in}\quad \Omega^c,
\end{array}
\end{equation}
where $p\in(0,p_\alpha^*)$ and $\delta_0$ is the Dirac mass at $0$. We recall the following result dealing with the convolution operator $\ast$ in Lorentz spaces $L^{p,q}(\R^N)$ (see \cite{O'N}).

\begin{proposition}\label{Conv}
Let $1\leq p_1, q_1, p_2, q_2\leq \infty$ and suppose $\frac{1}{p_1}+\frac{1}{p_2}>1$. If
$f\in L^{p_1,q_1}(\R^N)$ and $g\in L^{p_2,q_2}(\R^N)$, then $f\ast g\in L^{r,s}(\R^N)$ with $\frac{1}{r}=\frac{1}{p_1}+\frac{1}{p_2}-1$, $\frac{1}{q_1}+\frac{1}{q_2}\geq \frac{1}{s}$ and there holds
\begin{equation}\label{lorentz}
\norm {f\ast g}_{L^{r,s}(\R^N)}\leq 3r\norm f_{L^{p_1,q_1}(\R^N)}\norm g_{L^{p_2,q_2}(\R^N)}.
\end{equation}
In the particular case of Marcinkiewicz spaces $L^{p,\infty}(\R^N)=M^p(\R^N)$, the result takes the form
\begin{equation}\label{marcink}
\norm {f\ast g}_{M^{r}(\R^N)}\leq 3r\norm f_{M^{p_1}(\R^N)}\norm g_{M^{p_2}(\R^N)}.
\end{equation}
\end{proposition}

\begin{proposition}\label{pr 3}
Assume that $0<p<p^*_\alpha$ and $u$ is a nonnegative weak solution of (\ref{eq4.1}). Then
\begin{equation}\label{4.2.3}
0\le u\le \mathbb{G}_\alpha[\delta_0],
\end{equation}
$|\nabla u|\in L^{\infty}_{loc}(\Omega\setminus\{0\})$ and $u$ is a classical solution of
\begin{equation}\label{eq1.3}
 \arraycolsep=1pt
\begin{array}{lll}
 (-\Delta)^\alpha  u+ |\nabla u|^p=0\quad & \rm{in}\quad\Omega\setminus\{0\},\\[2mm]
 \phantom{   (-\Delta)^\alpha +|\nabla u|^p}
u=0\quad & \rm{in}\quad \Omega^c.
\end{array}
\end{equation}
\end{proposition}
{\it Proof.}  Since $0<p<p^*_\alpha$, (\ref{eq4.1}) admits a solution. Estimate (\ref{4.2.3}) is a particular case of  (\ref{1.5}). We pick a point $a\in \Omega\setminus\{0\}$ and
consider a finite sequence $\{r_j\}_{j=0}^\kappa$ such that $0<r_\kappa<r_{\kappa-1}<...<r_0$ and $\bar B_{r_0}(a)\subset \Omega\setminus\{0\}$. We set $d_j=r_{j-1}-r_{j}$, $j=1,...\kappa$. By (\ref{L9}) with $\beta=0$, it follows that
\begin{equation}\label{4.2.4}
\int_{\Omega}\left(u+|\nabla u|^p\right)dx\leq c_{20}.
\end{equation}
Let $\{\eta_n\}\subset\C^\infty_0(\R^N)$ be a sequence of radially decreasing and symmetric mollifiers such that supp$(\eta_n)\subset B_{\varepsilon_n}(0)$ and $\varepsilon_n\leq \frac12\min\{\rho(a)-r_0,|a|-r_0\}$ and $u_n=u\ast\eta_n$. Since
$$\eta_n\ast(-\Delta)^{\alpha}\xi= (-\Delta)^{\alpha}(\xi\ast\eta_n)
$$
by Fourier analysis and
$$\int_{\R^N}\!\left(u(-\Delta)^{\alpha}(\xi\ast\eta_n) \!+\xi\ast\eta_n|\nabla u|^p\right)\!dx=\!
\int_{\R^N}\!\left(u\ast\eta_n(-\Delta)^{\alpha}\xi+\eta_n\ast |\nabla u|^p\xi\right)\! dx
$$
because $\eta_n$ is radially  symmetric, it follows that $u_n$ is a classical solution of
\begin{equation}\label{eq1.4}
 \arraycolsep=1pt
\begin{array}{lll}
(-\Delta)^{\alpha}u_n+ |\nabla u|^p\ast\eta_n=\eta_n\qquad&\mbox{in }\ \Omega_n,\\
\phantom{(-\Delta)^{\alpha}+ |\nabla u|^p\ast\eta_n}u_n=0&\mbox{in }\ \Omega^c_n,
\end{array}
\end{equation}
where $\Omega_n=\{x\in \R^N:\mbox{dist} (x, \Omega)<\varepsilon_n\}$. We denote by ${G}_{\alpha,n}(x,y)$ the Green kernel of $(-\Delta)^{\alpha}$ in $\Omega_n$ and by $\mathbb{G}_{\alpha,n}$ the Green operator.
Set
$f_n=\eta_n- |\nabla u|^p\ast\eta_n$, then $u_n=\mathbb{G}_{\alpha,n}[f_n]$. If we set $f_{n,r_0}=f_n\chi_{B_{r_0}(a)}$, $\tilde f_{n,r_0}=f_n-f_{n,r_0}$, we have
$$\begin{array}{lll}
\D \partial_{x_i}u_n(x)=\int_{\Omega_n}\partial_{x_i}{G}_{\alpha,n}(x,y)f_n(y)dy
\\[4mm]\phantom{\partial_{x_i}u_n(x)}
\D =\int_{\Omega_n}\partial_{x_i}{G}_{\alpha,n}(x,y)f_{n,r_0}(y)dy+\int_{\Omega_n}\partial_{x_i}{G}_{\alpha,n}(x,y)\tilde f_{n,r_0}(y)dy
\\[4mm]\phantom{\partial_{x_i}u_n(x)}
=v_{n,r_0}(x)+\tilde v_{n,r_0}(x),
\end{array}$$
where
$$v_{n,r_0}(x)=\int_{B_{r_0}(a)}\partial_{x_i}{G}_{\alpha,n}(x,y)f_{n}(y)dy=-\int_{B_{r_0}(a)}\partial_{x_i}{G}_{\alpha,n}(x,y)|\nabla u|^p\ast\eta_n(y) dy$$
and
$$\tilde v_{n,r_0}(x)=\int_{\Omega_n\setminus B_{r_0}(a)}\partial_{x_i}{G}_{\alpha,n}(x,y)f_{n}(y)dy.
$$
We set $\rho_n(x)=\mbox{dist} (x,\Omega_n^c)$, then by (\ref{annex 01}) and (\ref{cap7 annex 01}), we have
$$|\partial_{x_i}{G}_{\alpha,n}(x,y)|\leq c_4N\max\left\{\frac1{|x-y|^{N-2\alpha+1}},
\frac{\rho^{-1}_n(x)}{|x-y|^{N-2\alpha}}\right\}.
$$
 Thus, if $x\in B_{r_1}(a)$ and $y\in \Omega_n\setminus B_{r_0}(a)$, then $\rho_n(x)>d_1$ and $|x-y|>d_1$,
\begin{equation}\label{V1}
|\tilde v_{n,r_0}(x)|\leq c_{21}\int_{\Omega_n\setminus B_{r_0}(a)}f_{n}(y)dy
\leq c_{20}c_{21},
\end{equation}
where $c_{21}>0$ depends  on $d_1^{-N+2\alpha-1}$, $N$ and $\alpha$.
Furthermore, if $x\in B_{r_1}(a)$ and $y\in B_{r_0}(a)$,
\begin{equation}\label{V}
|\partial_{x_i}{G}_{\alpha,n}(x,y)|\leq \frac{c_{4}N}{|x-y|^{N-2\alpha+1}}.
\end{equation}
We have already use the fact that $y\mapsto |y|^{2\alpha-N-1}\in L_{loc}^{q_1}(\R^N)$
with $q_1\in(\max\{1,p\},p_\alpha^*)$. Since $f_n$ is uniformly bounded in $L^1(\Omega)$, there exists
$c_{22}>0$ such that
\begin{equation}\label{V2}
\norm{v_{n,r_0}}_{M^{q_1}(B_{r_1}(a))}\leq c_{22}.
\end{equation}
 Combined with (\ref{V1}), it yields
\begin{equation}\label{V3}
\norm{|\nabla u|^p\ast\eta_n}_{M^{\frac{q_1}{p}}(B_{r_1}(a))}\leq c_{23}.
\end{equation}
Next we set $f_{n,r_1}=f_n\chi_{B_{r_1}(a)}$ and $\tilde f_{n,r_1}=f_n-f_{n,r_1}$. Then
$$\partial_{x_i}u_n=v_{n,r_1}+\tilde v_{n,r_1},
$$
where
$$v_{n,r_1}(x)=\int_{B_{r_1}(a)}\partial_{x_i}G_\alpha(x,y)f_n(y) dy=-\int_{B_{r_1}(a)}\partial_{x_i}G_\alpha(x,y)|\nabla u|^p\ast\eta_n(y) dy
$$
and
$$\tilde v_{n,r_1}(x)=\int_{\Omega_n\setminus B_{r_1}(a)}\partial_{x_i}G_\alpha(x,y)f_n(y) dy
$$
Clearly $\tilde v_{n,r_1}(x)$ is uniformly bounded in $B_{r_2}(a)$ by a constant $c_{24}$ depending on the structural constants and $d_2=r_1-r_2$. Estimate (\ref{V}) holds if we assume $x\in B_{r_2}(a)$ and $y\in B_{r_1}(a)$. Therefore,
$$|v_{n,r_1}(x)|\leq c_{4}N\int_{B_{r_1}(a)}\frac{|\nabla u|^p\ast\eta_n(y) }{|x-y|^{N-2\alpha+1}}dy.
$$
We derive from Proposition \ref{Conv}
$$\norm{v_{n,r_1}}_{M^{q_2}(B_{r_2}(a))}\leq c_{24}\norm{|\nabla u|^p\ast\eta_n}_{M^{\frac{q_1}{p}}(B_{r_1}(a))},
$$
with
\begin{equation}\label{V4}
\frac{1}{q_2}=\frac{p}{q_1}+\frac{1}{q_1}-1.
\end{equation}
Notice that $q_2>q_1$. Therefore
\begin{equation}\label{V5}
\norm{|\nabla u|^p\ast\eta_n}_{M^{\frac{q_2}{p}}(B_{r_2}(a))}\leq c_{25}.
\end{equation}
We iterate this construction and obtain the existence of constants $c_{j}$ such that
\begin{equation}\label{V4-j}
\norm{|\nabla u|^p\ast\eta_n}_{M^{\frac{q_j}{p}}(B_{r_j}(a))}\leq \bar c_{j},\qquad\forall j=1,2,....
\end{equation}
We pick $q_1=\frac12(p^*_\alpha+p)$ if $p>1$ or $q_1=\frac12(p^*_\alpha+1)$ if $p\in(0,1]$
\begin{equation}\label{V6}
\frac{1}{q_{j+1}}=\frac{p}{q_j}+\frac{1}{q_1}-1.
\end{equation}
If $p=1$,  there exists $j_0\in\N$ such that $q_{j_0}>0$ and $q_{j_0+1}\le0$.\\
If $p\in(0,p_\alpha^*)\setminus\{1\}$, let $\ell=\frac{q_1-1}{q_1(p-1)}$, then $\ell=p\ell+\frac{1}{q_1}-1$, thus
\begin{equation}\label{V7}\begin{array}{lll}\D
\frac{1}{q_{j+1}}=\ell +p^j\left(\frac{1}{q_{1}}-\ell\right)=\ell -p^j\frac{q_1-p}{q_1(p-1)}.
\end{array}\end{equation}
Therefore there exists $j_0$ such that $q_{j_0}>0$ and $q_{j_0+1}\leq 0$. This implies
\begin{equation}\label{V6}
\norm{|\nabla u|^p\ast\eta_n}_{L^s(B_{r_{j_0+1}}(a))}\leq c_{26},\qquad\forall s<\infty
\end{equation}
and
\begin{equation}\label{V7}
\norm{|\nabla u|^p\ast\eta_n}_{L^\infty(B_{r_{j_0+2}}(a))}\leq c_{27},
\end{equation}
with $c_{27}$ independent of $n$. Letting $n\to\infty$ infers
\begin{equation}\label{V8}
\norm{{\nabla u}}_{L^\infty(B_{r_{j_0+2}}(a))}\leq c^{\frac{1}{p}}_{27}.
\end{equation}
Combining this estimate with (\ref{4.2.3}) and using \cite[Corollary 2.5]{RS} which states
\begin{equation}\label{V9}\begin{array}{lll}
\norm{u}_{C^\beta(B_{r_{j_0+3}}(a))}\leq c\left(\norm u_{L^1(\R^N,\frac{dx}{1+|x|^{N+2\alpha}})}\right.\\[2mm]\phantom{----------}\left.+\norm u_{L^\infty(B_{r_{j_0+2}}(a))}+\norm{{\nabla u}}_{L^\infty(B_{r_{j_0+2}}(a))}\right),
\end{array}\end{equation}
for any $\beta<2\alpha$, we obtain that $u$ remains bounded in $C^{1+\varepsilon}(K)$ for any compact set
$K\subset \Omega\setminus\{0\}$ and some $\varepsilon>0$. Using now \cite[Corollary 2.4]{RS}, we obtain that
$C^{2\alpha+\varepsilon'}(\Omega\setminus\{0\})$ for $0<\varepsilon'<\varepsilon$. Futhermore $u$ is continuous up to $\partial\Omega$. As a consequence it is a strong solution in  $\Omega\setminus\{0\}$.\hfill$\Box$
\medskip



In the next result we give a pointwise estimate of $\nabla u$ for a positive solution $u$ of (\ref{eq4.1}).
\begin{proposition}\label{pr X} Assume that $R=\frac12dist(0,\partial\Omega)$, $p\in(0,p^*_\alpha)$ and $u$ is a nonnegative weak solution of (\ref{eq4.1}). Then there exists $c_{28}>0$ depending on $R$, $p$ and $\alpha$ such that
\begin{equation}\label{X1}
|\nabla u(x)|\leq c_{28}|x|^{2\alpha-N-1},\qquad\forall x\in \bar B_{R/4}(0)\setminus\{0\}.
\end{equation}
\end{proposition}
{\it Proof.}  Up to a change of variable we can assume that $R=1$. For $0<|x|\leq 1$, there exists $b\in (0,1)$ such that $b/2\leq |x|\leq b$. We set
$$u_b(y)=b^{N-2\alpha}u(by). $$
Then
$$(-\Delta)^{\alpha}u_b+b^{N+p(2\alpha-N-1)}|\nabla u_b|^p=0\qquad\mbox{in }\ \Omega_b:=b^{-1}\Omega.
$$
Using \cite[Corollary 2.5]{RS} with $\beta<2\alpha$, for any $a\in\Omega_b$ such that $|a|=3/4$, there holds
\begin{equation}\label{X2}\begin{array}{lll}\norm{u_b}_{C^\beta(B_{\frac{3}{16}}(a))}\leq c_{29}\left(\norm {u_b}_{
L^1(\R^N,\frac{dx}{1+|y|^{N+2\alpha}})}+\norm {u_b}_{L^\infty({B_{\frac{3}{8}}}(a))}\right.\\\phantom{----------}
\left.+b^{N+p(2\alpha-N-1)}\norm{|\nabla u_b|^p}_{L^\infty({B_{\frac{3}{8}}}(a))}^{\phantom{}}\right).
\end{array}\end{equation}
Furthermore, by the same argument as in Proposition \ref{pr 3},
\begin{equation}\label{X3}
\norm{|\nabla u_b|^p}_{L^\infty({B_{\frac{3}{8}}}(a))}\leq c_{30}\int_{\Omega_b}|\nabla u_b(y)|^pdy
=c_{30}b^{p(N+1-2\alpha)-N}\int_{\Omega}|\nabla u(x)|^pdx,
\end{equation}
and from (\ref{4.2.3}) and (\ref{annex 01})
$$u(x)\leq G_\alpha(x,0)\leq \frac{c_4}{|x|^{N-2\alpha}}\Longrightarrow u_b(y)\leq \frac{c_4}{|y|^{N-2\alpha}}.
$$
Then
$$\norm {u_b}_{L^1(\R^N,\frac{dy}{1+|y|^{N+2\alpha}})}\leq
c_4\int_{\R^N}\frac{dy}{|y|^{N-2\alpha}(1+|y|)^{N+2\alpha}}
=c_{31}.$$
If we take $\beta=1$, which is possible since $\alpha>1/2$, we derive
$$|\nabla u_b(a)|\leq c_{32}\Longrightarrow |\nabla u(ba)|\leq c^{-1}_{32}b^{2\alpha-N-1}
$$
In particular, with $|b|=4|x|/3$ we derive (\ref{X1}) with $c_{28}=c_{32}^{-1}(\frac{4}{3})^{2\alpha-N-1}$.\hfill$\Box$\medskip

We denote
\begin{equation}\label{X4}
c_{N,\alpha}=\lim_{x\to 0}|x|^{N-2\alpha}G_{\alpha}(x,0).
\end{equation}
It is well known that $c_{N,\alpha}>0$ does not depend on the domain $\Omega$ and, by the maximum principle,
$G_{\alpha}(x,0)\leq c_{N,\alpha}|x|^{2\alpha-N}$ in $\Omega\setminus\{0\}$.

\begin{teo}\label{teo 2}
Let $\Omega$ be an open bounded $C^2$ domain containing $0$, $\alpha\in(\frac12,1)$ and  $0<p<p^*_\alpha$.
If $u$ is a positive solution of problem (\ref{eq4.1}) and $\bar B_R(0)\subset \Omega$, it satisfies \smallskip

\noindent(i)  if $\frac{2\alpha}{N-2\alpha+1}<p<p^*_\alpha$,
$$0<\frac{ c_{N,\alpha}}{|x|^{N-2\alpha}}-u(x)\le
\frac{c_{33}}{|x|^{(N-2\alpha+1)p-2\alpha}},\quad x\in B_{R/4}(0)\setminus\{0\};
$$

\noindent(ii)  if $p=\frac{2\alpha}{N-2\alpha+1}$,
$$0<\frac{ c_{N,\alpha}}{|x|^{N-2\alpha}}-u(x)\le -c_{33}\ln (|x|),\quad x\in B_{R/4}(0)\setminus\{0\};$$

\noindent(iii) if $0<p<\frac{2\alpha}{N-2\alpha+1}$,
$$ 0<\frac{ c_{N,\alpha}}{|x|^{N-2\alpha}}-u(x)\le c_{33},\quad x\in B_{R/4}(0)\setminus\{0\},$$
where $c_{33}$ depends on $N$, $p$, $\alpha$ and $R$.
\\ Furthermore, if $1\leq p<p^*_\alpha$, this solution is unique.

\end{teo}

\noindent {\it Proof.} The existence of a nonnegative weak solution is a consequence of the subriticality assumption;  the fact that this solution is a classical solution in $\Omega\setminus\{0\}$ derives from Proposition \ref{pr 3}. It follows by (\ref{4.2.3})  and (\ref{4.2.4}) that for any $x\in \Omega\setminus\{0\}$,
\begin{equation}\label{X5}\begin{array}{ll}\displaystyle
\frac{ c_{N,\alpha}}{|x|^{N-2\alpha}}-u(x)\leq \int_{\Omega}G_{\alpha}(x,y)|\nabla u(y)|^pdy\\[4mm]
\displaystyle\phantom{\frac{ c_{N,\alpha}}{|x|^{N-2\alpha}}-u(x)}\leq c^p_{28}c_{4}\int_{B_{\frac{R}{4}}(0)}|x-y|^{2\alpha-N}|y|^{p(2\alpha-N-1)}dy +c_{34}\norm{\nabla u}_{L^{p}(\Omega)}\\[4mm]\displaystyle\phantom{\frac{ c_{N,\alpha}}{|x|^{N-2\alpha}}-u(x)}
\le c_{35}\left[\int_{B_{\frac{R}{4}}(0)}|x-y|^{2\alpha-N}|y|^{p(2\alpha-N-1)}dy +1\right]
\end{array}\end{equation}
where $c_{34},c_{35}>0$ depend on $N$, $p$ and $\alpha$. Next we assume $0<|x|\leq\frac{R}{16}$.\smallskip

\noindent {\it Case: $\frac{2\alpha}{N-2\alpha+1}<p<p^*_\alpha$}.
We can write
$$\int_{B_{\frac{R}{4}}(0)}|x-y|^{2\alpha-N}|y|^{p(2\alpha-N-1)}dy=E_1+E_2
$$
with
$$E_1=\int_{B_{\frac{R}{4}(0)}\setminus B_{\frac{R}{8}}(0)}|x-y|^{2\alpha-N}|y|^{p(2\alpha-N-1)}dy\leq c_{36},
$$
where $c_{36}>0$ depends on $N$, $\alpha$, $p$ and $R$ and
$$\begin{array}{lll}\displaystyle
E_2=\int_{B_{\frac{R}{8}}(0)}|x-y|^{2\alpha-N}|y|^{p(2\alpha-N-1)}dy\\[4mm]
\displaystyle\phantom{E_1}=|x|^{2\alpha-p(N+1-2\alpha)}
\int_{B_{\frac{R}{8|x|}}(0)}|\xi-\zeta|^{2\alpha-N}|\zeta|^{p(2\alpha-N-1)}d\zeta\\[4mm]
\displaystyle\phantom{E_1}
\leq \int_{|\zeta|>2}|\xi-\zeta|^{2\alpha-N}|\zeta|^{p(2\alpha-N-1)}d\zeta
\end{array}$$
with $\xi=x/|x|$. Since $2\alpha-N<0$, $|\xi-\zeta|^{2\alpha-N}\leq (|\zeta|-1)^{2\alpha-N}$, then
$$E_2\leq c_N\int_{2}^\infty (r-1)^{2\alpha-N} r^{p(2\alpha-N-1)+N-1}dr=c_{37}.
$$
Thus (i) follows.
\smallskip

\noindent {\it Case: $\frac{2\alpha}{N-2\alpha+1}=p$}. We see that
$$E_2=\int_{B_{\frac{R}{8|x|}}(0)}|\xi-\zeta|^{2\alpha-N}|\zeta|^{-2\alpha}d\zeta,
$$
then clearly
$$E_2=-\ln |x|+o(1)\quad\mbox{when }\ |x|\to 0.
$$
Thus (ii) follows.
\smallskip

\noindent {\it Case: $0<p<\frac{2\alpha}{N-2\alpha+1}$}. We have that
$$E_2=\int_{B_{\frac{R}{8|x|}}(0)}|\xi-\zeta|^{2\alpha-N}|\zeta|^{-2\alpha}d\zeta
=c_{29}|x|^{ p(N+1-2\alpha)-2\alpha}+o(1)\quad\mbox{when }\  |x|\to 0.
$$
Thus (iii) follows.

 Uniqueness in the case $1\leq p<p^*_\alpha$, is very standard, since if   $u_1$ and $u_2$ are two positive solutions of (\ref{eq4.1}), they satisfies
$$\lim_{x\to 0}\frac{u_1(x)}{u_2(x)}=1.
$$
Then, for any $\varepsilon>0$, $u_{1,\varepsilon}:=(1+\varepsilon)u_1$ is a supersolution which dominates $u_2$ near $0$, it follows by the maximum principle that $w:=u_2-(1+\varepsilon)u_1$ satisfies
$$(-\Delta)^\alpha w+|\nabla u_2|^p-|\nabla u_{1,\varepsilon}|^p\leq 0
$$
since $w$ is negative near $0$ and vanishes on $\partial\Omega$, if it is not always negative, there would exists $x_0\in \Omega\setminus\{0\}$ such that $w(x_0)$ reaches a maximum and $|\nabla u_2(x_0)|=|\nabla u_{1,\varepsilon}(x_0)|$, thus $(-\Delta)^\alpha w(x_0)\leq 0$, contradiction.
\hfill$\Box$\medskip

\noindent{\it Remark.} If $0<p<1$, the nonlinearity is not convex  and uniqueness does hold only if two solutions
$u_1$ and $u_2$ satisfy
$$\lim_{x\to 0}(u_1(x)-u_2(x))=0.
$$

\end{document}